\documentclass[a4paper,12pt,leqno]{amsart} 
\usepackage{amsmath,amsthm,amssymb}  
\usepackage{graphicx}

\numberwithin{equation}{section}
\numberwithin{figure}{section}

\theoremstyle{plain}
\theoremstyle{definition}  
\newcommand{\C}{\mathbb C}   
\newcommand{\R}{\mathbb R}
\newcommand{\Z}{\mathbb Z}

\newcommand{\al}{\alpha}
\newcommand{\be}{\beta} 
\newcommand{\ga}{\gamma}
\newcommand{\de}{\delta}
\newcommand{\la}{\lambda}
 
\newcommand{\La}{\Lambda}
\newcommand{\eps}{\epsilon}

\newcommand{\De}{\Delta}
\renewcommand{\th}{\theta}
\newcommand{\om}{\omega}

\newcommand{\Ga}{\Gamma}

\newcommand{\Si}{\Sigma}

\newcommand{\diag}{\text{diag}}
\newcommand{\Hol}{\text{Hol}}

\newcommand{\SL}{\textrm{SL}}
\renewcommand{\sl}{\frak s\frak l}

\newcommand{\GL}{\textrm{GL}}
\renewcommand{\O}{\textrm{O}}

\newcommand{\g}{{\frak g}}
\newcommand{\h}{{\frak h}}

\newcommand{\no}{\noindent}

\newcommand{\st}{\ \vert\ }   
\renewcommand{\ll}{\lq\lq}
\newcommand{\rr}{\rq\rq\ }
\newcommand{\rrr}{\rq\rq}  
\renewcommand{\b}{\partial}

\newcommand{\bp}{\begin{pmatrix}} 
\newcommand{\ep}{\end{pmatrix}} 
\newcommand{\bsp}{\left(\begin{smallmatrix}} 
\newcommand{\esp}{\end{smallmatrix}\right)}

\newcommand{\ttb}{ {t\bar t}  }

\renewcommand{\i}{ {\scriptscriptstyle\sqrt{-1}}\, }

\newcommand{\hb}{\hbar}

\newcommand{\Mzw}{ M^{(0)} }

\begin{document}     

\title[Quantum cohomology]{Quantum cohomology: is it still relevant?
}  
   
\author{Martin A. Guest}    

\date{}   

\maketitle 

\begin{abstract}
This article,  intended for a general mathematical audience, 
is an informal review of some of the many
interesting links which have developed between
quantum cohomology and \ll classical\rr mathematics. It is based on a
talk given at the Autumn Meeting of the Mathematical Society of Japan 
in September 2021.
\end{abstract}

 \section{Preamble}\label{1}
 
 If you are a geometer, you probably know that
 quantum cohomology came into prominence over 20 years ago, but you might not
 have followed developments since then.  If you are not a geometer, you
 might ask first \ll What is quantum cohomology?\rrr.  To answer this, this article will begin
 with some brief historical remarks. 
 
 Even before that, as preliminary orientation, I should say
 that quantum cohomology belongs to the \ll structure building\rr kind
 of mathematics, not to the \ll problem solving\rr kind. As far
 as I know, quantum cohomology was not motivated by the desire to solve 
 any particular problem.  There are no famous open problems in the subject. There are
 no easily-stated conjectures. 
 
 In fact there are relatively few general theorems (at least, in the traditional 20th century sense of the word).
 More precisely, most theorems in the subject apply to rather special situations, or even to special examples.  In this sense, 
 quantum cohomology theory shares some of the features of integrable systems theory (which
 will appear later in the article), in particular an emphasis on examples. Both could be
 described as unfinished theories --- they are \ll work in progress\rrr.
  
 On the other hand, quantum cohomology has revealed many remarkable
 links between previously-unrelated areas of mathematics. Many of these have led to very concrete results. Thus, quantum cohomology has
 provided answers to interesting problems --- but, in many cases, the problems were recognised only after the answers had been obtained. 
 We shall describe some of these links in order to demonstrate that
 quantum cohomology continues to have a strong impact --- that it is indeed
 still relevant.

\no{\em Acknowledgements:\  }  The author was partially supported by JSPS grant 18H03668.

 \section{Some history}\label{2}
  
  Here we sketch the general ideas (and some sociology) of
  quantum cohomology. Technical material is postponed to the
  following sections.
  
  The \ll quantum\rr in quantum cohomology comes from physics,
  of course. In quantum theory it is necessary to consider integrals over \ll the space of all paths\rrr, but a rigorous mathematical formulation of such integrals (over infinite-dimensional spaces) is a
  difficult problem. 
  
  Various approaches have been tried, with varying degrees of success.  
  One approach is to ignore the original problem, and focus on direct
  proofs of its expected consequences. 
  As a
  simple analogy, consider Cauchy's Theorem and the theory of residues in complex analysis. By integrating meromorphic functions with infinitely many poles, one can obtain interesting identities, such as
 \[
 \pi^2/6 = 1 + 1/2^2 + 1/3^2 + \cdots.
 \]
 A rigorous proof of Cauchy's Theorem requires knowledge of topology and function theory, 
 but direct proofs of identities --- such as this one --- can often be given without using complex analysis.
 
Another approach is to impose a (large) symmetry group and then hope for a localization principle which reduces the problem to the analysis of (small) fixed point sets.  
 
 From such considerations, physicists
 were led to consider the space
 \[
 \Hol(\Si,M) = \{ f:\Si\to M\st f \text{ is holomorphic} \}
 \]
of holomorphic maps from a Riemann surface $\Si$ to
a K\"ahler manifold $M$.  The
connected components of this space are generally noncompact, and may have singularities, but at least they are finite-dimensional.  

For example, holomorphic maps from the Riemann sphere $\Si=\C P^1$ to $n$-dimensional complex projective space $M=\C P^n$ are given by $(n+1)$-tuples of polynomials
\[
[p_0(z);\dots;p_n(z)]
\]
and the connected components of $ \Hol(\Si,M)$ are indexed by $d=0,1,2,\dots$ where 
$d=\text{max}\{\deg p_0,\dots,\deg p_n\}$
(and where $p_0,\dots,p_n$ are assumed to have no common factor).
It is plausible that a theory of integration over such spaces might exist.  

In fact, a mathematical foundation for such a theory of integration --- much more generally, for pseudo-holomorphic curves in symplectic manifolds --- was constructed by 
 Mikhael Gromov
in the 1980's. This involved a deep study of the (non)compactness of
the space.   Algebraic geometers provided another approach: in the 1990's 
 Maxim Kontsevich
introduced the concept of stable map, in order to construct an \ll improved\rr version of $\Hol(\Si,M)$.
Both approaches --- moduli spaces of curves, in symplectic geometry and algebraic geometry --- are nontrivial and became major areas of research. 
Some references for this are \cite{Ma99}, \cite{McSa04}.

\subsection{Phase I}

Quantum cohomology was born (into the mathematical world) in the 1990's as \ll intersection theory\rr on such moduli spaces.  The value of this kind of intersection theory had been demonstrated in the 1980's by 
 Simon Donaldson in the context of Yang-Mills theory on four-dimensional manifolds. Physicists had already carried out many ad hoc, but very promising,  calculations. Thus, the scene was set.

When the holomorphic maps are constant, quantum cohomology reduces to 
intersection theory on the target manifold $M$ itself, i.e.\ ordinary cohomology. Thus, 
quantum cohomology appeared to be a rather natural generalization.  On the other hand, 
quantum cohomology is very different from  ordinary cohomology in several respects.
First, it is not (in any naive sense) functorial. This means that
quantum cohomology is difficult to compute systematically.  Second, it does not (directly) measure any topological quantity, as the symplectic/complex structure plays an essential role.  
 
 In view of the lack of functoriality, early research on quantum cohomology concentrated on \ll Gromov-Witten invariants\rrr. These are the raw data of intersection theory.  They amount to giving a list of all possible intersections. This is not as bad as it sounds, as it suffices to \ll choose a basis\rr and take all possible intersections of the basis elements (a countable set). Moreover, for physicists, this is
 a perfectly satisfactory method --- it is what they always do.  There is
 also a mathematical antecedent: 
 computing Gromov-Witten invariants \ll one by one\rr is the analogue of
 Schubert calculus.  Schubert calculus could be described as \ll intersection theory on
 Grassmannian manifolds, before the invention of cohomology\rrr.  
 
 Incidentally, quantum
 Schubert calculus was one of the early applications of
 quantum cohomology.  
It has interesting applications to topics far from topology and physics, for example the \ll multiplicative eigenvalue problem\rr in linear algebra \cite{AgWo98}. (This is an example of discovering the answer first, then recognising the problem later.)

\subsection{Phase II}

The story so far could be called Phase I of quantum cohomology theory. We regard
Phase II as the next period of development, when these rather unruly 
Gromov-Witten invariants were organised more systematically.  This involved something quite new:  a relation with differential equations.  (The relation was discovered by
physicists, and in this sense pre-dates Phase I, but we are ordering the \ll Phases\rr 
in logical order, not time order.) 

The famous predictions of Mirror Symmetry --- counting rational curves in Calabi-Yau manifolds ---  came from solving differential equations.
This was carried out in a 1991 paper \cite{CdGP91} 
by 
 Philip Candelas, Xenia de la Ossa, Paul  Green, and   Linda Parkes.
But we shall take the starting point of Phase II as the address by
 Alexander Givental
at the 1994 ICM in Z\"urich \cite{Gi95}.  

Givental introduced the \ll quantum differential equations\rrr,
a linear system of partial differential equations.  
The key property
of the quantum differential equations is that there exists a solution
whose series expansion is a generating function of the (generalized) Gromov-Witten invariants.  

A simple example of this type of differential equation is
\begin{equation}\label{sto}
\tfrac{d^2u}{dx^2} + \tfrac1x \tfrac{du}{dx} - u = 0.
\end{equation}
This example was studied in a 19th century article \cite{St1857} by Stokes.
It is essentially the Bessel equation, but we mention it here because the Stokes Phenomenon for this
o.d.e.\  will enter the story later, in section \ref{4}.

A great advantage of this point of view is that the theory of  linear differential equations is more familiar than intersection theory on moduli spaces, and differential equations are closer to the original physics.  On the other hand it is hard to imagine which differential equations are related to quantum cohomology (clearly, not arbitrary ones), and it
is hard to imagine how to extract  intersection-theoretic information. Nevertheless, Mirror Symmetry gave us
some nontrivial examples.  Exploring this phenomenon was the task of Phase II.  

The dependent and independent variables of the quantum differential equations were coordinates
on cohomology vector spaces. This means cohomology with complex coefficients;  the coefficients of 
quantum differential equations were holomorphic.  To a topologist, doing calculus on cohomology spaces is
unusual, to say the least.  Givental called this phenomenon \ll homological geometry\rrr.

In this new world of homological geometry, the foundations of  quantum cohomology (Phase I),
at least 
for \ll nice\rr target spaces $M$, could be taken for granted.  Relying on these hard-won foundations, armed only with
undergraduate knowledge of differential equations (and a little algebraic topology), many more researchers could
plunge into quantum cohomology theory.

This was a very fruitful exercise, because quantum differential equations could
be postulated --- that is, guessed --- even for target spaces $M$ where the
foundations were not yet secure.  For example: target spaces with singularities,
target spaces with boundary, noncompact target spaces, infinite-dimensional
target spaces, and so on. 
Even if the differential equations did not quite fit the
(known) geometry, one could take the self-righteous view that
the differential equations were \ll correct\rr and the geometry \ll wrong\rrr. 
In turn, these efforts led to further
development of the foundations. 

It was a major achievement to give a mathematical treatment of the physicists' 
original Mirror Symmetry predictions for (certain) Calabi-Yau
manifolds. This was done by Givental \cite{Gi96} and also
by 
 Bong Lian, Kefeng Liu, and  Shing-Tung Yau
\cite{LLY1}.
The work of Lian-Liu-Yau also led to
intriguing connections with arithmetic and modular forms.

Toric varieties (and their subvarieties) provided the most
important class of target spaces for these results. 
The lack of functoriality of quantum cohomology made it necessary to focus on
examples, and toric varieties came to the rescue here.  

Contrarily, and perhaps with \ll tongue in cheek\rrr,  it could be said that Mirror Symmetry
came to the rescue of toric varieties.  Certainly the theory of toric varieties
became much more widely known, because of Mirror Symmetry.
It was a happy coincidence, and another unexpected benefit of quantum cohomology.

A good reference for Phase II of quantum cohomology theory is \cite{CoKa99}.

Before going on to Phase III, it should be mentioned that the symplectic approach saw its own dramatic developments. The key words here were Floer theory and Witten's version of Morse theory. 
Floer theory for the loop space $\La M$ is closely related to the quantum cohomology of $M$,
but Floer theory developed in its own way.  In particular the case of a target manifold
$M$ with Lagrangian submanifold $L$, and Gromov-Witten invariants based on
maps from $(D,\b D)$ to $(M,L)$, became a major theme,
developed by
 Kenji Fukaya, Yong-Geun Oh, Hiroshi Ohta, and  Kaoru Ono \cite{FOOO09}.
Here $D$ is a two-dimensional
disk with boundary $\b D$.

\subsection{Phase III}

In this (admittedly very subjective) narrative, Phase III started with
another ICM talk, namely that of 
 Boris Dubrovin
in 1998 \cite{Du98}.  Here another
mathematical leap forward was made: from 
{\em linear} to {\em nonlinear} differential equations.  

As always in this story, physicists had anticipated such a leap, with
their WDVV (Witten-Dijkgraaf-Verlinde-Verlinde) equations
and tt* (topological-antitopological fusion) equations.  
These equations involve more general
kinds of quantum cohomology:  

--- WDVV describes \ll big\rr quantum cohomology

--- tt* describes \ll quantum cohomology with a real structure\rrr.

Both are complicated and difficult to appreciate, partly because
of the physicists' preference for coordinates and tensor notation.  Dubrovin
introduced a more intrinsic point of view for both, based on his
concept of Frobenius manifold \cite{Du96}.

In the context of quantum cohomology, the Frobenius manifold is the 
cohomology vector space of the target manifold, endowed with
further differential geometric data. (After all, having started to do calculus on 
cohomology vector spaces, why not introduce connections and curvature as well?) 

Independently from  quantum cohomology, this kind of structure had already been investigated by
 Kyoji Saito
in the 1980's, in the context of singularity theory. Thus
Frobenius manifolds linked quantum cohomology with singularity theory in a systematic way.
To mathematicians, this was remarkable, but it was entirely consistent
with the physical principles of Mirror Symmetry.

The stage was set for another unexpected application of quantum cohomology:  the mere
existence of a target space $M$ (or rather, its quantum cohomology) implies
the existence of a solution of a highly nontrivial nonlinear p.d.e.  

For example, Dubrovin showed that the WDVV equation for $\C P^2$ reduces to a case of the Painlev\'e VI equation, so 
the \ll big\rr quantum cohomology of 
$\C P^2$ corresponds to a certain solution of this equation. While this was by no means the first time that Painlev\'e equations had
appeared in geometry and physics, it was a radically new example.
Dubrovin also began the study of integrable hierarchies generated by the quantum differential equations,
in analogy with the KdV hierarchy which is generated by the Schr\"odinger equation.
Such hierarchies represent a further generalization, which could be called
\ll very big\rr  quantum cohomology.

\subsection{Next steps}

Many of the research areas started by, or influenced by, quantum cohomology during
the past 20 years have reached a certain level of maturity; they have become part
of mainstream mathematics.  As often happens, some have become
highly specialized.   As often happens, specialists in one area 
may not be able to talk to specialists in another.

A glance at the lists of ICM talks shows
that quantum cohomology was most prominent in the period 1990-2002. That
was 20 years ago. The question in the title
\ll Is it still relevant?\rr
was suggested by this.

\section{Some terminology}\label{3}

We shall now give some simple examples related to Phase II and Phase III, 
to prepare some terminology for the next section.  We use only undergraduate-level material in this section.

The simplest version of the quantum differential equation for $M=\C P^n$
(and  $\Si=\C P^1$)
 is
\begin{equation}\label{cpn}
(\hb\b)^{n+1} y = qy
\end{equation}
where $y=y(q)$, $q$ is a complex variable, $\b=q\frac{d}{dq}$, and $\hb$ is
a parameter.  The case $n=1$ is Stokes' equation (\ref{sto}) above, if we put $y=u$, $q=x^2$ (and $\hb=2$).

This is related to the cohomology of $\C P^n$ in the following way.
The cohomology vector space $H^\ast(\C P^n;\C)$ has
$1,b,b^2,\dots,b^n$ as a basis, where $b$ is an (additive)
generator of $H^2(\C P^n;\C)$.  It is generated
multiplicatively by $b$, with the relation $b^{n+1}=0$. 
The quantum  cohomology is equal to $H^\ast(\C P^n;\C)$ as
a vector space, but it has a \ll deformed\rr product operation,
leading to the relation $b^{n+1}=q$.  (Ordinary cohomology
is recovered from quantum  cohomology by setting $q=0$.) The relation
$b^{n+1}-q$ corresponds to
the operator  $(\hb\b)^{n+1} - q$
in an obvious way. 

Although $q=0$ is a singular point for the quantum differential equation,
a basis $y^{[0]},\dots,y^{[n]}$ of solutions (near $q=0$) may easily be found by the
Frobenius Method.  Givental observed that, if we write
\[
J=(y^{[0]},\dots,y^{[n]}) \ \longleftrightarrow\ 
y^{[0]}1 + y^{[1]}b +    \cdots + y^{[n]} b^n,
\]
then this Frobenius solution is given by the attractive formula
\[
J(q)=
q^{ b/\hb  } 
\ 
\sum_{k=0}^\infty 
\frac{q^k}
{[
(b+\hb)(b+2\hb)\cdots(b+k\hb)
]}.
\]
At first this looks strange, as $b$ is a cohomology class, but it makes sense
as $b^{n+1}=0$ (in cohomology).  With this
understanding, it is easy to verify that
$
(\hb\b)^{n+1} J = qJ,
$
as asserted.  There are good reasons to identify
\[
q
\ \longleftrightarrow\ 
qb \in H^2(\C P^n;\C),
\]
so we have a map $J:H^2(\C P^n;\C)\to H^\ast(\C P^n;\C)$
(strictly speaking, a multivalued map, as $q^{ b/\hb  } = e^{ (b/\hb) \log q  }$ gives rise to
logarithms).  The solution of the quantum differential equation is
a cohomology-valued function on a cohomology vector space!

For $M=\C P^n$ the relation between the quantum cohomology ring (i.e.\ the 
cohomology vector space with the deformed product) and the 
quantum differential equation is deceptively simple.  In particular we do not really need $\hb$ there. A less obvious example is the case of
a (smooth) hypersurface $M$ of degree $3$ in $\C P^4$.  There is a
generator $b\in H^2(M;\C)\cong \C$ which satisfies
$b^4=27qb^2$, but the quantum differential equation is
\[
\left(
(\hb\b)^4 - 27q (\hb\b)^2 - 27 \hb q (\hb\b) - 6 \hb^2 q
\right)
y=0.
\]
Replacing $\hb\b$ by $b$ in the differential operator gives $b^4-27qb^2-27\hb q b -6 \hb^2 q$, and this gives $b^4-27qb^2$
{\em after putting $\hb=0$.}  Here the inverse procedure is not immediately visible.  However, it turns out that there is such a procedure (see \cite{Gu05}), and here the parameter $\hb$ plays an essential role.  

More generally, if $\dim H^2(M;\C)=r$, then we have a system
of linear p.d.e.\ in variables $q_1,\dots,q_r$.  A new source
of difficulty appears:  while an o.d.e.\ $Py=0$ of order $k$ always has a (local) solution
space of dimension $k$,  there is no such dimension formula in terms of the orders $k_1,k_2,\dots$ of
the operators $P_1,P_2,\dots$ of a system.  In the situation of the
quantum differential equations this dimension must be equal to $\dim H^\ast(M;\C)$;
knowledge of this dimension is a nontrivial property of quantum cohomology.

In fact,  systems of linear p.d.e.\ with (nonzero, but) finite-dimensional
local solution space are quite rare.  They correspond to flat
connections in vector bundles, or D-modules of finite rank. Examples
of this correspondence (related to geometry in general, and quantum
cohomology in particular) are explored 
in detail in \cite{Gu08}.

Recall that a connection in a vector bundle on $(q_1,\dots,q_r)$-space
may be expressed locally as $\nabla=d+\om$, where $\om$ is a matrix-valued $1$-form. It is
said to be flat if $d\om+\om\wedge\om=0$. In the case of
quantum cohomology (or a Frobenius manifold), the connection
is called the Givental connection (or Dubrovin connection).

 In general, if the coefficients of
a connection form are written in terms of functions $f_1(q_1,\dots,q_r),f_2(q_1,\dots,q_r),\dots$
then the flatness condition is a system of nonlinear p.d.e.\ for  $f_1,f_2,\dots$
and there is a close relation between this nonlinear p.d.e.\ and the original
linear p.d.e.   
The theory of \ll integrable p.d.e.\rr amounts to studying this kind of relation.
For example, the (nonlinear) KdV equation is related
to the (linear) Schr\"odinger equation in exactly this way.

Quantum cohomology (in the sense of most of this article, i.e.\ small quantum cohomology of Fano manifolds) fits into this nonlinear p.d.e.\  framework, but only in a trivial way. This
is because the entries of the Dubrovin connection form are just polynomials on
$H^2(M;\C)$. For small quantum cohomology, the WDVV equations are trivial.
However, big quantum cohomology involves functions on 
$H^\ast(M;\C)$, which are \ll known\rr  only to the extent that they are solutions of the WDVV equations. 
Later we shall encounter another nonlinear system, the tt* equations, which is
nontrivial even for small quantum cohomology.

\section{Some recent developments}\label{4}

After all this preparation we can review some topics of current research (restricted by the interests and knowledge of the author).  

In the past 10 years, Fano manifolds (like $\C P^n$) have played
a greater role than Calabi-Yau manifolds, partly because Calabi-Yau manifolds were
the main focus in the early years, and partly because Fano manifolds require new methods. While Calabi-Yau manifolds are closely related to variations of Hodge structure (VHS), Fano manifolds require \ll semi-infinite\rr
or \ll nonabelian\rr Hodge structures ($\frac\infty 2$-VHS).
In terms of the quantum differential equations, the new feature of
Fano manifolds is that irregular singularities appear.  

The Hodge-theoretic aspects are 
a huge challenge, and major achievements have been made by 
 Claude Sabbah,
and by
 Takuro Mochizuki,
and their collaborators.  

The irregular singularity aspects
also suggest new directions of research.  This holds even
in the simplest situation where the quantum differential equation is an o.d.e.\ --- although
the classical theory goes back to the 19th century, the relation with quantum cohomology is relatively unexplored.  We shall focus on this aspect for the rest of the article.

\subsection{Stokes data of the  quantum differential equation}\label{4.1}

The quantum differential equation (\ref{cpn}) for $M=\C P^n$ has a regular
singularity at $q=0$ and an irregular
singularity at $q=\infty$. Near $q=0$ we have seen a nice series expansion, related to Gromov-Witten invariants. Near $q=\infty$ we certainly have
(local, possibly multivalued) solutions, but \ll series expansions at infinity\rr
always diverge. This is the Stokes Phenomenon: such series
represent only {\em asymptotic expansions} of solutions, and only on specific {\em Stokes sectors} of the complex plane.

Remarkably (and non-intuitively) the asymptotic expansion --- which depends only on the equation --- determines a unique solution {\em on a given Stokes sector.} 
Two such solutions on overlapping sectors must differ by a constant matrix; this
is called a Stokes matrix. Analytically continuing a given solution one circuit around $q=\infty$, we
pass through a finite number of sectors, so we
see that the monodromy of that solution is (essentially) a product of Stokes matrices.  

In general, Stokes matrices are difficult to compute, and their entries are usually highly transcendental, but in the
case of the quantum differential equations we might
expect them to contain geometric/physical information. This is indeed the case. For $M=\C P^n$
the Stokes data reduces to $n$ real numbers $s_1,\dots,s_n$.  These were computed by Dubrovin for $n=2$, and by 
 Davide Guzzetti
for general $n$, and the answer is surprisingly simple: $s_k=\binom{n+1}{k}$ (see \cite{Du98}). Physicists had also arrived at this fact by \ll counting solitons\rr in the supersymmetric sigma model of $\C P^n$.
Clearly something interesting was happening.

A further piece of information can be extracted from the differential equation,
by comparing the Frobenius solution near $q=0$ with (any of) the Stokes solutions near $q=\infty$. 
In the differential equations literature,
the matrices which occur here are called {\em connection matrices.}  For
$M=\C P^n$, Dubrovin computed these as well (see \cite{Du98}). 

A remarkable
geometrical interpretation of this was found by
 Hiroshi Iritani
(see \cite{Ir09}, \cite{Ir16}),
which conjecturally also applies to more general Fano manifolds $M$.
It involves the
\ll gamma class\rr
\[
\hat\Ga_M = \prod_{i=1}^n \,\Ga(1+x_i)
\]
where $x_1,\dots,x_n$ are the Chern roots of the
complex tangent bundle $TM$, and where $\Ga(1+x_i)$
is interpreted as the (Taylor expansion of) the gamma function.
This is very much in the spirit of Givental's homological geometry.
Moreover, it  strengthens a surprising link
with \ll helices\rr of holomorphic vector bundles on $M$ (see \cite{Za96}).
This has been developed further in \cite{GGI16}, \cite{CDGXX} --- another major new direction of research based on quantum cohomology.

\subsection{Isomonodromy}\label{4.2}

The Stokes data of the quantum differential equation can be approached most directly by converting the
quantum differential equation (a p.d.e.\ in $q_1,\dots,q_r$ with parameter $\hb$)
to an o.d.e.\ in $\hb$ with parameters $q_1,\dots,q_r$.  This \ll trick\rr is made possible
by homogeneity. 
An important property of the o.d.e.\ in $\hb$ is that it is isomonodromic: its 
monodromy data (Stokes and connection matrices) do not depend on $q_1,\dots,q_r$.
(The reverse is also true, i.e.\  the monodromy data of the
p.d.e.\ in $q_1,\dots,q_r$ is independent of $\hb$, but this is obvious.)

For big quantum cohomology the isomonodromy property shows that the monodromy data gives
{\em conserved quantities} of solutions of the WDVV equations. This provides another link with the theory of integrable systems. Indeed, Dubrovin proposed this as
a general approach to the study of Frobenius manifolds. 
There is a well-developed method --- the Riemann Hilbert Method --- for
studying the relation between  monodromy data of linear o.d.e.\
and solutions of nonlinear p.d.e.\  (in the case of
the Painlev\'e equations, a good reference is \cite{FIKN06}).

As we have mentioned, the WDVV equation for $\C P^2$ reduces to a case of the Painlev\'e VI equation, so this example fits well.  
Nevertheless,  the general WDVV equations are a \ll can of worms\rrr.  
On the other hand, if we restrict to small quantum cohomology, the WDVV equations become trivial. 
Fortunately, between these two extremes, there is a situation of intermediate difficulty, which appears when we consider Frobenius manifolds with \ll real structure\rrr.  The nonlinear p.d.e.\ here is the (system of) tt* equations.  

\subsection{The topological-antitopological fusion equations.}\label{4.3}

The tt* equations were introduced by
 Sergio Cecotti
and
 Cumrun Vafa
in the context of supersymmetric quantum field theory \cite{CeVa91}, 
\cite{CeVa92}. Dubrovin formulated these equations as an isomonodromic system
\cite{Du93}, so the Riemann-Hilbert Method of  \cite{FIKN06} can, in principle, be applied here.

Before stating the equations, which may look unmotivated, we can state
the underlying (mathematical) idea, which is simple.  A Frobenius manifold is a holomorphic object, but, if it has
a real structure, then there is a
\ll complex conjugate\rr Frobenius manifold, which is an antiholomorphic object.  
Now, a Frobenius manifold has, as part of its definition, a \ll holomorphic metric\rrr.  Thus,
the real structure produces a Hermitian metric.  The tt* equations
are the equations for this metric --- the tt* metric.

For Frobenius manifolds, e.g.\ big quantum cohomology, the
tt* equations are extremely complicated (they are worse than
the WDVV equations, which they contain).  
However, even after restricting to small quantum cohomology (where
the WDVV equations are trivial), some nonlinear equations remain. As
Dubrovin pointed out in \cite{Du93}, these nonlinear
equations are very familiar to differential geometers: they
are the equations for pluriharmonic maps from $(q_1,\dots,q_r)$-space
to the Riemannian symmetric space $\GL_{n+1}\C/\O_{n+1}$
(where $n+1$ is the dimension of the Frobenius manifold).
In fact, the pluriharmonic map represents the
variation of Hodge structure which was mentioned earlier.

When $r=1$, a pluriharmonic map
is just a harmonic map. Here the theory of harmonic maps from
surfaces into symmetric spaces (which was comprehensively
developed by differential geometers in the 1980's) provides an effective tool.  We shall
focus on this situation.

The tt* equations are then
\begin{equation}\label{ceva}
\tfrac{\b}{\b \bar t}
\left (
g \tfrac{\b}{\b t} g^{-1}
\right) 
-  [C,g C^\dagger g^{-1}]=0,
\end{equation}
where $C$ is the (holomorphic) chiral matrix of the theory, $C^\dagger$ is
the conjugate-transpose of $C$, and $g^{-1}$ is the Hermitian matrix representing the tt* metric.   
In the situation of quantum cohomology, $C$ is the matrix of quantum multiplication by
a generator $b\in H^2(M;\C)$ (in the situation of singularity
theory, $C$ is the analogous matrix of multiplication in the Milnor ring). 

Let us now specialize to the case $M=\C P^n$. Then, with respect to suitable bases,
we have

\begin{equation*}
C=
\left(
\begin{array}{c|c|c|c}
 & & & \ q
 \\
\hline
1\ &  &   &  
\\
\hline
  & \ddots  &  & 
\\
\hline
 & &\ 1\ & 
\end{array}
\right),
\end{equation*}
and the tt* metric can
be written
\[
g=\diag(e^{-2w_0},\dots,e^{-2w_n}),
\]
where $w_0,\dots,w_n$ are real-valued functions of $q$, and where
$q=\al t^\be$ for suitable constants $\al,\be$.
The homogeneity
property of quantum cohomology implies that $w_i$ depends only on $\vert q\vert$. 
Then 
the tt* equations 
become 
\begin{equation}\label{ttt}
 2(w_i)_{\ttb}=-e^{2(w_{i+1}-w_{i})} + e^{2(w_{i}-w_{i-1})}, \ 
 i=0,1,\dots,n
\end{equation}
together with the extra condition $w_i+w_{n-i}=0$.   This is
a version of the periodic Toda equations. (We
interpret $w_{n+1},w_{-1}$ here as $w_0,w_n$ respectively.)

These \ll Toda type\rr tt* equations were, in fact, one of
the main examples considered by Cecotti and Vafa (formula (7.4) in \cite{CeVa91}). We
call them the tt*-Toda equations.

Physically, a solution is a massive deformation of a conformal field theory, and the existence of such a deformation says something about that theory.
Cecotti and Vafa made a series of conjectures about the solutions:
 
 ---there should exist (globally smooth) solutions $w$ on $\C^\ast$
 which are radial, i.e.\ $w=w(\vert t\vert)$ 

---these solutions should be characterized by asymptotic data at $t=0$
(the \ll ultra-violet point\rrr; here the data is equivalent to the chiral charges, essentially the holomorphic
matrix function $C$)

---these solutions should equally be characterized by asymptotic data at $t=\infty$,
(the \ll infra-red point\rrr; here the data is equivalent to the Stokes parameters $s_i$,
or soliton multiplicities, about which we shall say more later on).

\no A global solution can then be interpreted as a \ll flow\rr between these two kinds of data (the renormalization group flow).

The Toda equation is a very familiar p.d.e., which has been well studied 
in
the theory of integrable systems (as well as in differential geometry). 
In addition to the methods of these areas, we
have available the Riemann-Hilbert Method from isomonodromy theory.
Unfortunately none of these methods work directly,  because it
is an important physical requirement that $w_i$ should be defined for all $t\ne 0$ in the complex plane. Obtaining globally defined solutions is the main obstacle.  

Nevertheless, when Cecotti
and Vafa introduced these equations, some information was
available for $n=1$, as this case had been studied in detail
by Barry McCoy, Craig Tracy, and  Tai Tsun Wu
in their pioneering work \cite{MTW77} on the Ising Model.  
For $n=1$ equation (\ref{ttt}) reduces to the
sinh-Gordon equation 
\begin{equation}\label{sG}
(w_0)_{\ttb}=\sinh 4w_0.
\end{equation} 
Moreover, because
of the  radial condition, it is equivalent to a special
case of the 
Painlev\'e III equation. From their work it was known
that radial solutions $w_0$ on $\C^\ast$ are in
one to one correspondence with
real numbers $\ga_0\in[-1,1]$, where
$2w_0(\vert t\vert) \sim \ga_0 \log \vert t\vert$ as $t\to 0$.
From isomonodromy theory it was known that the
Stokes parameter $s_1$ is given by
\begin{equation}\label{sGsto}
s_1= 2\sin \frac \pi2 \ga_0 \ \in\  [-2,2]
\end{equation}
(as $n=1$ there
is only one Stokes parameter). 

At this point we should clarify what we mean by \ll the tt* equations for
$\C P^n$\rrr.  In fact (the quantum cohomology of) $\C P^n$ corresponds to just one solution of (\ref{ttt});
the other solutions are also expected to have
geometric/physical significance, but this would
not necessarily be the quantum cohomology of anything.

In the case $n=1$, the solution corresponding to the quantum cohomology of $\C P^1$
is given by $\ga_0=1$ (or $\ga_0=-1$). Although quantum cohomology was new
in physics at that time (and still unknown in mathematics), 
Cecotti and Vafa were aware of this fact.  On the basis of this, and a small
number of similar examples, and much physical reasoning, they proposed that only the solutions with
{\em integer Stokes data} represent \ll realistic\rr physical models. 
Cecotti and Vafa used this hypothesis to propose a classification of 
supersymmetric field theories \cite{CeVa92}. 
From (\ref{sGsto}) we can observe that, in the case $n=1$, the Stokes parameter $s_1$ is an integer only for $\ga_0=0,\pm \frac13, \pm 1$. 
The solution with $\ga_0=\pm \frac13$ corresponds to the Ising Model,
and the solution with $\ga_0=0$ corresponds to a trivial model (here $w_0\equiv 0$). 
Thus, all solutions with $s_1\in\Z$ are meaningful.

The above predictions were investigated for the tt*-Toda equation (\ref{ttt})  in the
series of papers
\cite{GuLi12}, 
\cite{GuLi14}, \cite{GIL1}, \cite{GIL2}, \cite{GIL3}, \cite{MoXX}, \cite{Mo14}.
A summary of the main results can be found in section 3 of \cite{Gu21}. Amongst these are:

(i) solutions of (\ref{ttt}) on $\C^\ast$ are in
one to one correspondence with $(n+1)$-tuples
$(\ga_0,\dots,\ga_n)$ such that $\ga_{i+1}-\ga_i\ge -2$ (and $\ga_i+\ga_{n-i}=0$);
one has
$2w_i(\vert t\vert) \sim \ga_i \log \vert t\vert$ as $t\to 0$

(ii) the corresponding Stokes parameters $s_1,\dots,s_n$ are 
the $i$-th symmetric functions of 
$
e^{ (n-\ga_0)\frac{\pi\i}{n+1}  },
e^{ (n-2-\ga_1)\frac{\pi\i}{n+1}  },
\dots,
e^{ (-n-\ga_n)\frac{\pi\i}{n+1}  }
$
(thus $s_i=s_{n-i+1}$).

\no The relation between (i) and the chiral data, and
between (ii) and the soliton data, will be explained later. We just
remark here that the asymptotics at $t=\infty$ is related to
the Stokes parameters $s_i$ by the formula
\begin{equation}\label{linearcomb}
-\tfrac 4{n+1}
\sum_{p=0}^{[\frac12(n-1)]} w_p(\vert t\vert) \sin \tfrac{(2p+1)k\pi}{n+1} \sim
 s_k\, F(L_k \vert t\vert),
 \quad
 t\to\infty
\end{equation}
where 
$
F(x)=\tfrac12(\pi x)^{-\frac12}e^{-2x}
$,
$
L_k=2\sin \tfrac k{n+1}\pi
$,
and 
$[\tfrac12(n+1)]$ means $\tfrac12(n+1)$ if $n$ is odd, $\tfrac12n$ if $n$ is even.

The solution corresponding to the quantum cohomology of $\C P^n$ is given by
$(\ga_0,\dots,\ga_n)= (n,n-2,\dots,-n)$.  Then all exponentials
in (ii) above are equal to $1$, so $s_k=\binom{n+1}{k}$,
in agreement with section \ref{4.1}.

The methods used to obtain the above results were conventional ones: harmonic map theory (specifically, loop group factorizations), isomonodromy theory, and p.d.e.\ theory.  Rather less expected was the fact that all three methods played essential roles --- if used in isolation, none of the methods would have been sufficient to give the complete picture. 

Although the global solutions were the main physical motivation, we note that these methods 
give information also on wider classes of solutions, in particular solutions defined \ll near zero\rr (i.e.\ on regions
of the form $0<\vert t\vert< \eps$), and to solutions defined \ll near infinity\rr (i.e.\ on regions
of the form $R<\vert t\vert< \infty$). It is of mathematical interest to consider
such solutions from the Hamiltonian/symplectic point of view,
as a possible generalization of the rich theory which exists in the case $n=1$ for
the Painlev\'e III equation (and for any Painlev\'e equation).  Some initial
results in this direction can be found in \cite{OdXX}.

\subsection{Lie-theoretic aspects}\label{4.4}

The Lie group $\SL_{n+1}\C$ has been our silent companion so far --- the space
$\C P^n$ is a homogeneous space of $\SL_{n+1}\C$, and equations
(\ref{ttt}) are (a real form of) the $\sl_{n+1}\C$-Toda equations. It
is natural to consider other Lie groups. Although this
Lie-theoretic direction has not been developed very extensively, it is already
clear that there are benefits in doing so, for physics
as well as mathematics.  Regarding quantum cohomology, for example, it
was shown in \cite{KaXX} that minuscule flag manifolds for $G$ play the
same role as $\C P^n$ for $\SL_{n+1}\C$.

The most comprehensive results (at present) concern the Lie-theoretic generalization
of the Stokes data.
In classical o.d.e.\ theory the Lie group $\SL_{n+1}\C$ is implicit, but the Stokes Phenomenon for other Lie groups $G$ was investigated by
 Philip Boalch
in his work \cite{Bo02} on hyper-K\"ahler moduli spaces of meromorphic $G$-connections. We make use of this in the next section.

\section{Some applications}\label{5}

We have not yet explained how the chiral matrix $C$ (and the Stokes data of the quantum differential equations)  are related to the tt* metric $g=e^{-2w}$.  Nor is this explained in the original work of
Cecotti and Vafa, because their Stokes data is for the isomonodromic version of the tt*-Toda equations,
not the quantum differential equations.  In fact the two kinds of Stokes data {\em coincide}, as far as the global solutions are concerned.  This was proved in \cite{GIL3}, as a consequence of the
differential geometric approach to the tt*-Toda equations, and it is this approach
which gives the most direct relation between $C$ and $g$. 

In the case of the tt*-Toda equations, the procedure is as follows.  We consider the (more general)
holomorphic matrix function
\begin{equation*}
C=
\left(
\begin{array}{c|c|c|c}
\vphantom{\dfrac12}
 & & &  q^{k_0}
 \\
\hline
\vphantom{\dfrac12}
q^{k_1}  &  &   &  \\
\hline
\vphantom{\dfrac12}
  & \  \ddots\  &  & 
\\
\hline
\vphantom{\dfrac12}
 & &q^{k_n}  &  \!
\end{array}
\right)
\end{equation*}
and then introduce a holomorphic $1$-form
\[
\om =\tfrac1\la C dq,\quad \la\in S^1.
\]
Regarding this as a  $1$-form with values in the loop algebra $\La\sl_{n+1}\C$, we can write
(locally) 
$\om=L^{-1}dL$, where $L$ is an  $\La\SL_{n+1}\C$-valued holomorphic function (of $q$ and $\la$).
For a suitable choice of $L$, and a suitable loop group Iwasawa factorization $L=L_\R L_+$, it
can be shown that the $1$-form $(L_\R)^{-1}dL_\R$ is exactly the $1$-form which
appears in the  isomonodromic formulation of the tt*-Toda equations.  The details (given in \cite{GIL3}) are
somewhat subtle, but the method itself is well known to differential geometers --- it is the DPW construction, or generalized Weierstrass representation, of the corresponding harmonic map. This
method can be used to produce local solutions of the tt*-Toda equations near $t=0$ (although significantly more work is needed to investigate when these are actually global solutions).

The global solutions of tt*-Toda equations were described in section \ref{4.3} in terms
of (i) asymptotic data, and (ii) Stokes data.  The chiral data $C$ gives a third description:

(iii) Fix $N>0$.  Solutions of (\ref{ttt}) on $\C^\ast$ are in
one to one correspondence with $(n+1)$-tuples
$(k_0,\dots,k_n)$ such that $k_i\ge-1$, 
$n+1+\sum_{i=0}^n k_i=N$
(and $k_i=k_{n-i+1}$ for $i=1,\dots,n$).

The variable $q$ is related to 
the variable $t$ of 
the tt*-Toda equations
by
$t=\tfrac{n+1}N q^{  \frac N{n+1} }$. The $(k_0,\dots,k_n)$ are related
to the $(\ga_0,\dots,\ga_n)$ of (i) by $\ga_i-\ga_{i-1}+1=\frac{n+1}N(k_i+1)$.
Thus the chiral data (represented by the $k_i$) is close to the asymptotic data at $t=0$ (represented by the $\ga_i$), as stated earlier.

From now on it will be convenient to introduce the notation
\[
m_i=-\tfrac12 \ga_i
\]
(thus we have $w_i\sim -m_i \log \vert t\vert$ as $t\to 0$).
The relation between the $k_i$ and the $m_i$ is given by 
\begin{equation}\label{mandk}
m_{i-1}-m_i+1=\tfrac{n+1}{N}(k_i+1).
\end{equation}

\subsection{The Coxeter Plane}\label{5.1}

We begin with a Lie-theoretic description of the Stokes data, 
taken from \cite{GH2}. We include it in this section as it could be
described as an application of the tt*-Toda equations to Lie theory. It
will also be fundamental for the applications to physics given in the next two subsections.

Let $\mathfrak g$ be a complex simple Lie algebra, with corresponding simply-connected Lie group $G$.
Let $\al_1,\dots,\al_l\in\h^\ast$ be a choice of simple roots of $\g$ with respect to the Cartan subalgebra $\h$.  The Weyl group $W$ is the finite group generated 
by the reflections $r_\al$ in all root planes $\ker\al$, $\al\in\De$.  
The Coxeter element is the element $\ga=r_{\al_1}\dots r_{\al_l}$ of $W$.
 Its order is called the Coxeter number of $\g$, and we denote it by $s$.  
 We shall mainly be concerned with the case $\g=\sl_{n+1}\C$. Here
 $l=n$, $s=n+1$, and $W$ is the permutation group
on $n+1$ objects.
\begin{figure}[h]\label{sl5c}
\begin{center}
\vspace{-2cm}
\includegraphics[scale=0.28, trim= 0 200 0 0]{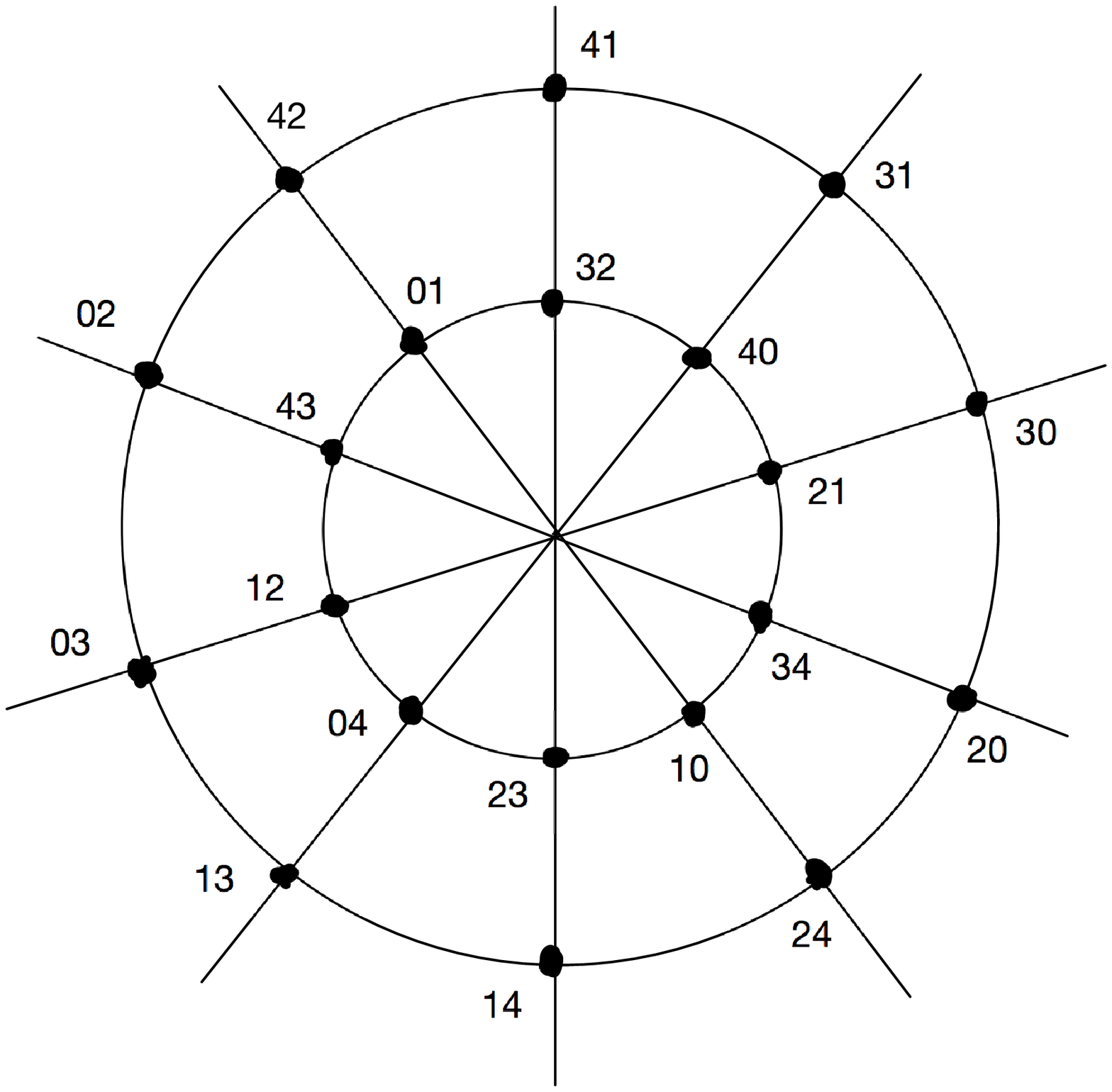}
\end{center}
\caption{The Coxeter Plane for $\g=\sl_5\C$.}
\end{figure}

The {\em Coxeter Plane} is the result of projecting $\De$  orthogonally onto a certain real plane in $\h^\ast$. 
Several versions of this definition can be found in the Appendix B of \cite{GH2}.
As an example, the Coxeter Plane for $\g=\sl_5\C$ is shown in Figure \ref{sl5c}.
Here there are $20$ roots $x_i-x_j$, $0\le i\ne j\le 4$, denoted by $ij$
in Figure \ref{sl5c}. The Coxeter element acts on the roots by
the permutation $(43210)$; there are $l=4$ orbits, each containing $s=5$ elements.

It was proved by Kostant (see \cite{Ko59}) that the
Coxeter element $\ga$ always acts on the set of roots $\De$ with $l$ orbits, each containing $s$ elements.  The Coxeter Plane provides a visualization of this fact.  Remarkably, it also 
provides a visualization of the Stokes data for the tt*-Toda equations.  Namely (see \cite{GH2}, \cite{GuHo21}):

 (a) The Coxeter Plane is a diagram of the Stokes sectors for the tt*-Toda equations.
 
(b) The Stokes matrices can be computed Lie-theoretically in terms of 
 a Lie group element 
 \[
 \Mzw=C(s_1,\dots,s_l)\in  \SL_{n+1}\C
 \]
where $C$ is
a  \ll  Steinberg cross-section\rr of the regular elements of $ \SL_{n+1}\C$.

To be more precise, the rays in the Coxeter Plane are the mid-points of the intersections of consecutive Stokes sectors, so for each ray we have a Stokes matrix.  To each ray in the Coxeter Plane is associated a set of roots,
and these roots determine the shape of the Stokes matrix. The Stokes matrices are related to the above matrix $\Mzw$ by a simple explicit formula (which we omit here). 

The corresponding asymptotic data at $t=0$ also has a Lie-theoretic interpretation.
Recall that the
solutions are parametrized by $(n+1)$-tuples 
$(m_0,\dots,m_n)$ satisfying

(*) \  $m_{i-1}-m_i+1\ge 0$  \ \  (equivalently, $k_i\ge-1$)

(**)\  $m_i+m_{n-i}=0$   \ \ (equivalently, $k_i=k_{n-i+1}$)

\no where 
$w_i\sim -m_i \log \vert t\vert$ as $t\to 0$. The
inequalities (*) define a convex polytope.
Let us write
\[
m=\diag(m_0, m_1,\dots,m_n),
\quad
\rho=\diag(\tfrac n2, \tfrac n2-1,\dots,-\tfrac n2).
\]
With this notation the semisimple part of
the matrix $\Mzw$ is
\[
\exp\, \tfrac{2\pi\i}{n+1}(m+\rho),
\]
and the convex polytope given by the points
$
\tfrac{1}{n+1}(m+\rho)
$
which satisfy (*) is the Fundamental Weyl Alcove of the Lie algebra.

We have already remarked that the 
solution corresponding to the quantum cohomology of $\C P^n$ is given by
$(\ga_0,\dots,\ga_n)= (n,n-2,\dots,-n)$; this corresponds to $m=-\rho$, which
is just the origin of the Fundamental Weyl Alcove.

\subsection{Particles and polytopes}

Using the material in the previous subsection, we can
show how the Coxeter Plane
and the tt*-Toda equations give a mathematical
foundation for certain field theory models proposed
by mathematical physicists in the 1990's (\cite{Fr91}, \cite{Do91}, \cite{Do92}).

These authors proposed (amongst other things) the
correspondence
\begin{align*}
\text{particle} \  &\leftrightarrow\   \text{ orbit of root in Coxeter Plane}
\\
\text{mass of particle} \  &\leftrightarrow\   \text{ distance of root from origin}
\end{align*}
(if $\g=\sl_{n+1}\C$ the mass of the particle corresponding
to the orbit of the root $x_i-x_j$ is $2\sin\vert i-j\vert \tfrac\pi{n+1}$).

They checked that these proposals (as well as the other things) were consistent with the
expected properties of a field theory. 

A variant of this proposal  was made (see \cite{FLMW91},  \cite{LeWa91})
for \ll polytopic models\rrr.  In this kind of model, a finite-dimensional
representation $\th$ of the Lie algebra $\g$ on a vector space $V$ is chosen, and the 
\ll polytope\rr is the polytope in $\h^\ast$ spanned by the weights of the representation.
The weight vectors (in $V$) are taken to be the vacua of the theory. 
In this theory, \ll solitonic particles\rr tunnel between vacua:   a soliton connects two vacua $v_i,v_j$ if
and only if the corresponding weights $\la_i,\la_j$
differ by a single root, i.e.\ $\la_i-\la_j\in\De$.
The physical characteristics of this particle are those of the root (in the model described above).

The discussion so far is purely algebraic (there is no differential equation).
However, the polytopic models include certain Landau-Ginzburg models. The
quantum cohomology of $\C P^n$ is of this type, with:
$
\th=\la_{n+1}
$
(standard representation of $\sl_{n+1}\C$). Thus one can expect a role for solitons 
in the quantum cohomology of $\C P^n$.

\begin{figure}[h]\label{sl4c}
\begin{center}
\vspace{-2.5cm}
\includegraphics[angle=90,origin=c,scale=0.3, trim= 500 50 -80 50]{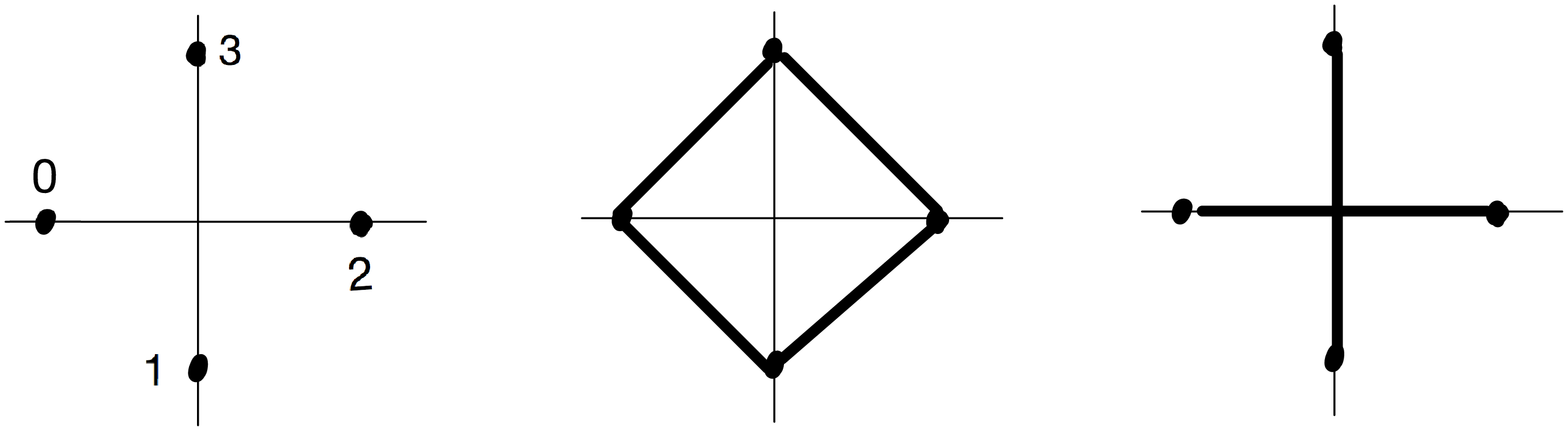}
\end{center}
\caption{Solitons for $\g=\sl_4\C$, $\th=\la_4$.}
\end{figure}

For example, in Figure \ref{sl4c}, the solitons are illustrated for $\sl_{4}\C$ with $\th=\la_4$.
The first part
shows the projections of the weights $x_0,x_1,x_2,x_3$, denoted by $0,1,2,3$ in
Figure \ref{sl4c}. The second part
shows (as heavy lines) the four solitons of type $x_0-x_1$ (and of mass $2\sin\frac \pi 4 = \sqrt 2$).
The third part shows the two solitons of  type $x_0-x_2$ (and of mass $2\sin\frac \pi 2 = 2$).
In this example, any two vacua are connected by a soliton.  
Other examples can be found in \cite{GuHo21} and in the original articles.

The \ll soliton particle interpretation\rr fits well with the tt*-Toda equations.  From a
physical point of view this was already implicit in the articles of Cecotti and Vafa. From a
mathematical point of view it arises from the formula (\ref{linearcomb}) 
on the asymptotics at $t=\infty$.
Namely, the linear combination on the {\em left hand side} of 
\[
-\tfrac 4{n+1}
\sum_{p=0}^{[\frac12(n-1)]} w_p \sin \tfrac{(2p+1)k\pi}{n+1} \sim
 s_k\, F(L_k \vert t\vert)
\]
corresponds to a certain basis vector of $\h$ (or $\h^\ast$) associated
to an orbit of the Coxeter group (see section 4.3 of \cite{Gu21}). 
Thus we can say that  the Stokes parameter $s_k$ on the
{\em right hand side}  is naturally associated 
to the $k$-th orbit, or particle. Physicists call $s_k$ the soliton multiplicity.

For any $n$, the model with $\th=\la_{n+1}$ is related to the
quantum cohomology of $\C P^n$.  If we choose
$
\th=\wedge^k \la_{n+1}
$
we obtain a different model, related to the
quantum cohomology of the Grassmannian $Gr_k(\C^{n+1})$.
Cecotti and Vafa used this to give a physical argument for an
\ll equivalence\rr
\[
\wedge^k QH^\ast( \C P^n ) \approx
QH^\ast( Gr_k(\C^{n+1}) )
\]
(more precisely, an equivalence of underlying field theories).
This tt* argument was explained in some detail in \cite{Bo95a}.

Later, mathematicians gave proofs of mathematically precise, but physically somewhat
artificial, interpretations of this isomorphism (which they regarded
as a special case of
the quantum Satake isomorphism, or abelian-nonabelian correspondence
(cf.\  \cite{GoMaXX}).

Our Lie-theoretic description of the solutions
of the tt*-Toda equations supports the
original physics argument, because
\begin{align*}
\text{solution with $m=-\rho$}  \ & 
\overset{\ \ \th=\la_{n+1}}{\longleftrightarrow}
\  \ QH^\ast( \C P^n )
\\
\text{solution with $m=-\rho$}   \ & 
\overset{ \ \ \ \th=\wedge^k\la_{n+1}}{\longleftrightarrow}
  QH^\ast( Gr_k(\C^{n+1}) )
\end{align*}
i.e.\  both $QH^\ast( \C P^n )$ and $ QH^\ast( Gr_k(\C^{n+1}) )$ arise from
the same solution of the tt*-Toda equations.

The Stokes matrices of the respective quantum differential equations
are different (they can be read off from $\Mzw$ and $\wedge^k\Mzw$
respectively), but the Stokes parameters $s_k=\binom{n+1}{k}$ are
the  same for $QH^\ast( \C P^n )$ and $ QH^\ast( Gr_k(\C^{n+1}) )$.
For further explanation we refer to \cite{Gu21}.

\subsection{Minimal models}

Recall that solutions of the tt*-Toda equations can be interpreted as harmonic maps;
in fact they are (rather special) examples of harmonic bundles. From this
point of view the chiral data
\begin{equation*}
\eta=
\left(
\begin{array}{c|c|c|c}
\vphantom{\dfrac12}
 & & &  z^{k_0}
 \\
\hline
\vphantom{\dfrac12}
z^{k_1}  &  &   &  \\
\hline
\vphantom{\dfrac12}
  & \  \ddots\  &  & 
\\
\hline
\vphantom{\dfrac12}
 & &z^{k_n}  &  \!
\end{array}
\right)
\end{equation*}
(more precisely the $1$-form $\eta dz$) can be interpreted as a Higgs field.  
Here we have 
$k_i\ge -1$,
$n+1+\sum_{i=0}^n k_i=N$,
$k_i=k_{n-i+1}$ for $i=1,\dots,n$.  As stated at the beginning of section \ref{5}, for fixed $N>0$,
such Higgs fields
parametrize 
solutions of the tt*-Toda equations.  

In this section, following Fredrickson and Neitzke \cite{FrNeXX},   we 
consider $\eta dz$ with 
$k_i\in \Z_{\ge 0}$,
$n+1+\sum_{i=0}^n k_i=N$,
such that $N$ is coprime to $k=\sum_{i=0}^n k_i$.  (We write $z$ instead of $q$, and $\eta$ instead of $C$, in recognition of the fact that we are moving away from the tt* interpretation. 
However, it should be noted that those Higgs fields with the additional
property 
$k_i=k_{n-i+1}$ for $i=1,\dots,n$ do form
 a dense subset of the set of all solutions of the tt*-Toda equations.)
 
Given (coprime) $k$ and $N$, there are only a finite number of Higgs fields of the above type. 
The authors of \cite{FrNeXX} describe these as those points of a certain moduli space
of Higgs field which are fixed under a
$\C^\ast$-action.  
This moduli space is associated to a certain four-dimensional supersymmetric quantum field theory, namely  Argyres-Douglas theory of type 
$(A_{n}, A_{k-1})$.

They observe
a \ll curious 1-1 correspondence\rr between these fixed points and certain representations of the
vertex algebra $W_{n+1}$.
The representations constitute the $(n+1,N)$ $W_{n+1}$ minimal model, a type of conformal field theory.
The
vertex algebra $W_{n+1}$ is known (or conjectured) to appear also in 
Argyres-Douglas theory of type 
$(A_{n}, A_{k-1})$, Fredrickson and Neitzke propose this as
a basis for relating the two models.

As an application of the Lie-theoretic Stokes formula
\[
\Mzw=C(s_1,\dots,s_l) \in  \SL_{n+1}\C
\]
we can show (see \cite{GuOtXX}) that there is
a rather direct mathematical path from $\eta dz$ to the representation.  

This depends on the theory of positive energy representations of the
affine Kac-Moody algebra
$\widehat{\sl}_{n+1}\C$. Irreducible representations of this type are parametrized
by pairs $(\La,l)$, where $l\in \Z_{\ge 0}$ and $\La$ is a dominant weight
of $\sl_{n+1}\C$ of level $l$. 

Let us review this notation briefly (cf.\ Remark 4.1 of \cite{GuOtXX}).  
First we denote the simple roots $x_0-x_1,\dots,x_{n-1}-x_n$
by $\al_1,\dots,\al_n$, then define the basic weights $\eps_1,\dots,\eps_n$ by the
condition $\langle \al_i,\eps_j\rangle = \de_{i,j}$.  The dominants weights are then given
by $\La=\sum_{i=1}^{k} k_i\eps_i$ where $k_i\in \Z_{\ge 0}$, and $\La$ is said
to have level $l$ if $\sum_{i=1}^{k} k_i\le l$.  

Let $P_+$ denote the set of dominants weights, and $P_l$ those of level $l$.
It is well known that
\[
P_l + \rho = P_+ \cap (l+n+1) \mathring A
\]
where $\mathring A$ denotes the interior of the
Fundamental Weyl Alcove
$A$.

The assumption $k_i\in \Z_{\ge 0}$ implies (see \cite{GH2}) (a) that
$\Mzw$ is semisimple, so it is conjugate to the diagonal matrix
\[
\exp\, {\tfrac{2\pi\i}{n+1}(m+\rho)},
\]
and also (b) that $\frac{1}{n+1}(m+\rho)$ is in $\mathring A$.

Now, it is easy to verify that our formula (\ref{mandk}) is equivalent to the formula
\begin{equation}\label{mandk2}
\tfrac N{n+1}(m+\rho) = \rho + \sum_{i=1}^n k_i\eps_i.
\end{equation}
We have $\sum_{i=1}^n k_i\le k=\sum_{i=0}^n k_i$. 
Thus, from the Stokes data $\Mzw$ we obtain the 
positive energy
representation corresponding to the pair $( \sum_{i=1}^n  k_i\eps_i,k)$.

It is well known (see \cite{BoSc93}, \cite{GuOtXX})
that $W_{n+1}$  intertwines with any such
representation, and that the effective central charge is then given by the formula
\[
\textstyle
c_{\text{eff}}=
n  -  12 \tfrac {n+1}N \,  \vert \sum_{i=1}^n k_i\eps_i - \tfrac{k}{n+1} \rho \vert^2.
\]
By (\ref{mandk2})  we have 
$ \sum_{i=1}^n k_i\eps_i - \tfrac{k}{n+1}\rho = \tfrac N{n+1}m $,
so
\[
c_{\text{eff}}
=n  -  12 \tfrac N{n+1} \,  \vert m \vert^2.
\]
This is the formula used by Fredrickson and Neitzke to relate 
representations of  $W_{n+1}$ and Higgs fields.

{\em

\noindent
Department of Mathematics\newline
Faculty of Science and Engineering\newline
Waseda University\newline
3-4-1 Okubo, Shinjuku, Tokyo 169-8555\newline
JAPAN
}
  
\end{document}